\documentclass[12pt,reqno]{amsart}
\usepackage{hyperref}
\usepackage{amsmath}
 \usepackage{bm,amsthm,amsfonts,latexsym,amssymb,commath} 

\usepackage{graphicx,float,multicol}

\usepackage{siunitx,caption}

\usepackage{threeparttable}[caption]
\usepackage{array, booktabs}
\begin{document}

\title[Modeling the endovenous laser ablation]{Analytical solutions in the modeling of the endovenous laser ablation}
\author{Luisa Consiglieri}
\address{Luisa Consiglieri, Independent Researcher Professor, European Union}
\urladdr{\href{http://sites.google.com/site/luisaconsiglieri}{http://sites.google.com/site/luisaconsiglieri}}

\begin{abstract} 
We model the operative treatment of incompetent truncal veins using endovenous laser ablation (EVLA).
Three differential equations, namely the diffusion,  the heat and the bioheat equations, are considered in the endovenous-perivenous multidomain,
describing the lumen, the vein wall, the tissue pad and the skin.
Exact solutions are provided. Our main concern is to accurate the heat source by taking the Beer--Lambert law into account in the irradiance
of the incident beam. To accurate the heat transfer at the skin boundary,  the Newton law of cooling is considered as a Robin boundary condition.
Open problems are presented.
\end{abstract}


\maketitle

\section{Introduction}
\label{intro}

Conventional  high ligation and stripping 
(crossectomy of the saphenofemoral junction (SFJ) with great saphenous vein (GSV) stripping)
 and radio-frequency (RF)  ablation therapies are passing their legacy to the new technologies such as the laser 
(and the ultrasound-guided foam sclerotherapy)
 in the treatment of varicose veins \cite{carr,etlik,rasm,theiva}.
Also the extremely high degree of recurrence occurring after previous ligation and stripping of the great saphenous vein
ask for the
safety and efficacy of endovenous laser ablation as a posterior treatment
for recurrent symptomatic saphenous insufficiency \cite{anch}.
In the recent years, several clinical follow-up studies
 (see, for instance, \cite{doganci,gale,gold,min,oh,pugg}) have reporting
different laser systems to treat incompetent GSV.
Indeed, the advantages of the endovenous laser surgery are more relevant than its complications \cite{ash}.
Clinical studies address failure as patency or recanalization of the GSV or residual
symptoms \cite{firo}, the occlusion, ulcer healing,
paresthesia rates and postoperative pain \cite{ozkan,rath}, as other adverse side effects \cite{shar,hode}.
Although it is unknown its cause even an ischemic
stroke following endovenous laser treatment of varicose veins is reported \cite{cagg}.

The fiber type is the single most significant factor related to treatment outcome \cite{prince}.
For effective endovenous laser (EVL)  therapy, the laser wavelength varies from
\(\lambda = \SI{810}{\nano\metre}\) \cite{shar,theiva} until to  \SI{980}{\nano\metre} \cite{schm},
at power settings of \SIrange{9}{17}{\watt} \cite{rath}
according to the diameter and severity of varicose veins. 
The 1470-nm wavelength EVL system successfully closes saphenous veins  
but not acts as "anesthesia-less thermal ablation technology" \cite{almeida}.
In \cite{proe}, the authors found that after EVLA with low-energy density, worse results and more relapses could be expected than
with higher-energy doses.

A vast literature has been playing a prominent and broad-spectrum role in the study of the temperature dependence
of the thermophysical and mechanical properties of biological tissues (see \cite{saccomandi} and the references therein).
The aim of some studies is to determine how the heat sink effect of the blood flow inside the
vessel may either be measured \cite{lc2012,cs,kotte},
distort coagulation volume during thermal therapies \cite{ho-shih,vuyl}, or  protect the
vessel wall in the proximity of an RF-assisted resection device \cite{g-suarez}. In this last work,
the tissue vaporization was modeled  by the enthalpy method, while
parametric studies were conducted in \cite{lc2012} to prove the blood flow has a cooling effect
during RF ablation treatment.
The present goal goes somewhat on the opposite direction: to study how the laser
wavelength behaves and the damage of blood vessels is influenced.

The finding of mathematical models  is essential to control the temperature, and to prevent postoperative complications. 
It includes handling
of mathematical problems generated by the application of models to real thermal conditions.
The combination Mordon's optical-thermal model with the presence
of a strongly absorbing carbonized blood layer on the fiber tip
is introduced in \cite{ruij} by neglecting the Arrhenius damage integral, and after  is developed in \cite{polue}.
We refer to \cite{marqa} a finite element modeling of the influence of air cooling that simulates the perisaphenous
subcutaneous tumescent saline solution infiltration.
Some analytical solutions to the bioheat transfer problem are studied  in  multiregion 
\cite{lc2013,lc2015}, where the Joule effect is assumed to be constant.

Here, we follow the whole path from the phenomenological interpretation of thermal therapy to scientific computing for producing simulations. 
Bearing this in mind, one mathematical model is stated and analytically solved such that the final method will be very fast to execute.
To validate our model, we compare our results  with experimental measurements.

\section{Mathematical model}
\label{math}

\begin{figure}
\begin{multicols}{2}
\centering
\includegraphics[width=0.5\textwidth]{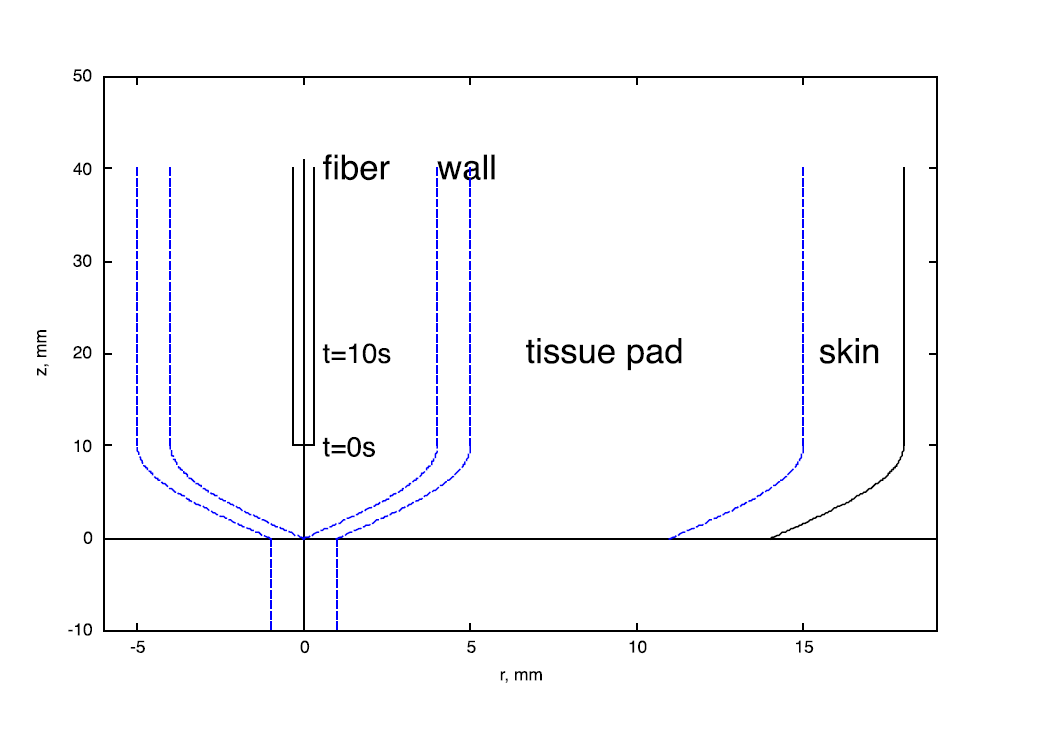}\\
\includegraphics[width=0.5\textwidth]{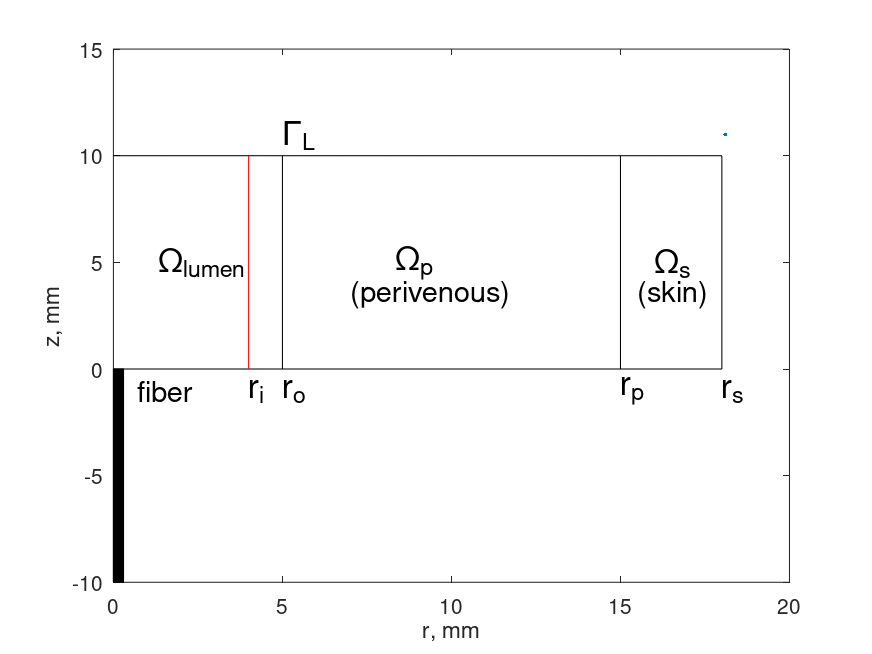}
\end{multicols}
\caption{\footnotesize Left: Schematic sagittal representation of the endovenous fiber at the initial instant of time \(t=0\).
Right: Schematic cylindrical  representations in 2D of the lumen
\( \Omega_\mathrm{lumen} =\{(x,y):\,x^2+y^2<r_{\mathrm{i}}^2 \} \times ]0;L[ \),
 the wall  \( \Omega_\mathrm{w} =\{(x,y):\,r_\mathrm{i}^2 <x^2+y^2<(r_\mathrm{i}+\varepsilon )^2 \}\times ]0;L[ \),
 the tissue pad \( \Omega_\mathrm{pad}= \{(x,y):\, (r_\mathrm{i}+\varepsilon )^2<x^2+y^2<r_\mathrm{p} ^2 \}\times ]0;L[\), 
and the skin  \( \Omega_\mathrm{skin}= \{(x,y):\, r_\mathrm{p}^2<x^2+y^2<r_\mathrm{s} ^2 \}\times ]0;L[\). }
\vspace{0.2cm}
\hrule
\label{vein}
\end{figure}

The optical laser fiber is inserted into the sheath so that the fiber tip extends
\SI{2}{\centi\metre} beyond the end of the sheath
to avoid the melting of the sheath material \cite{shar,schm}.
The geometry of the fiber--tissue system is assumed to be as follows (cf. Fig. \ref{vein}).
The fiber probe (with radius \(r_\mathrm{f}\))  is assumed to be centered in the middle of the vein,
where the vein segment \(\Omega_\mathrm{v}=\Omega_\mathrm{lumen}\cup \Omega_\mathrm{w}\) is assumed to be cylindrical with
\(r_\mathrm{i}\) and  \(\varepsilon\) being the inner  radius and the thickness of the venous wall,
respectively.
In vivo, the thickness of the venous wall is approximately one tenth the width of its blood column:
\(\varepsilon = r_\mathrm{i} /5\).

The saphenous vein is considered to be parallel to the skin surface  \cite{marqa}.
Then the complete domain \(\Omega = \Omega_\mathrm{v}\cup \Omega_\mathrm{pad}\cup \Omega_\mathrm{skin}\)
  is assumed to be constituted by axially half-cylindrical volumes, namely vein, perivenous and skin tissues.
The tissue around the venous may be considered homogenous. The thickness
of the perivenous and skin tissues are \(l_\mathrm{pad}=\SI{10}{\milli\metre}\) and \(l_\mathrm{skin}=\SI{3}{\milli\metre}\), respectively.

The procedure begins by inserting the laser sheath and positioning the bare tip below SFJ.
At time \(t=0\), the fiber tip is located at axial coordinate \(z_0\). 
The fiber is kept in a piecewise pull-back model  until the level of the knee.
A  pulling back of the laser fiber and the introduced catheter is  at a constant velocity \(v\), when the irradiation is activated,  and then they
may be pulled back by increments, during the off periods. The cycle is repeated until a desired fixed distance.
This  piecewise pull-back in the application of laser light is preferred because
 the manually made continuous pull-back 
 depends on the operating surgeon technique and experience. 
For each fixed  irradiation time \(t_\mathrm{end}>0\),
  \(L>0\) stands for the corresponding length of the treated vein segment.

 We assume that outer surface \(\Gamma_\mathrm{L}= \{(x,y):\ x^2+y^2< r_\mathrm{s}^2 \}\times \{\mathrm{L}\}\)
  is an insulating boundary, that is, there is no outflow.

\subsection{Diffusion approximation of the radiative transfer equation}

The light power emitted out of the fiber tip into the ambient
blood is scattered towards the vein wall and the surrounding tissue.
The fluence rate \(\phi\) [\si{\watt\per\square\metre}]  verifies the diffusion equation  (see, for instance, \cite{w-vang6})
\begin{equation}\label{light}
\frac{1}{\nu}\frac{\partial\phi}{\partial t}-D\Delta \phi +\mu_\mathrm{a}\phi=S
\qquad \mathrm{in }\ \Omega\times ]0 ; t_\mathrm{end} [,
\end{equation}
where:
\begin{description}
\item[ \(\nu=c/n\)] is the speed of light in the tissue [\si{\metre\per\second}], as determined by the relative refractive index \(n\);
\item[ \(D=\frac{1}{3(\mu_\mathrm{a}+\mu'_\mathrm{s})}\)] is the diffusion coefficient;
\item[ \(\mu_\mathrm{a}\)] is the absorption coefficient [\si{\per\metre}];
\item[ \(\mu'_\mathrm{s}=(1-g)\mu_\mathrm{s}\)] is the reduced scattering coefficient [\si{\per\metre}], with
\(g\) being the scattering anisotropy coefficient, and
\(\mu_\mathrm{s}\) being the scattering coefficient;
\item[ \(S\)] is the source of scattered photons [\si{\watt\per\cubic\metre}],
which represents the power injected in the unit volume.
\end{description}

Usually \(S\) is assumed to be a source point in order to use the Green functions in the determination of the
solution \(\phi\) of the diffusion equation \eqref{light}.
For a single spherically symmetric point source emitting \(P_\mathrm{laser}\) [\si{\watt}], the
form of \(\phi\) at a distance \(r\) from the source is \(P_\mathrm{laser}\exp(-\mu_\mathrm{eff}r) /(4\pi D r)\) \cite{mord},
where \(\mu_\mathrm{eff}=\sqrt{ 3 \mu_\mathrm{a} (\mu_\mathrm{a}+\mu'_\mathrm{s}) }\) represents the effective attenuation coefficient. 

Laser energy is delivered along the vein  with a continuous emission \cite{anch}, and a fiber pull-back is kept at about 1-cm increment every 10 seconds \cite{ash}. 
In this work, we consider 
\begin{equation}\label{source}
S (r,z,t)\equiv S(r_\mathrm{f},z,t)= \frac{P_\mathrm{laser}}{\pi r_\mathrm{f}^2 }  \frac{\mu_\mathrm{s}(\mu_\mathrm{t}+g\mu_\mathrm{a})}{\mu_\mathrm{a}+ \mu_\mathrm{s} '}
\exp[-\mu_\mathrm{t}(z+vt)] ,
\end{equation}
for all \((r,z,t)\in  [0; r_\mathrm{f}[\times ]- vt;L[\times [0;t_\mathrm{end}[\),
where   \(\mu_\mathrm{t}=\mu_\mathrm{a}+\mu_\mathrm{s}\) is the attenuation coefficient \textit{i.e.}  the reciprocal of the average
distance light travels before being scattered or absorbed by the medium.
The above  expression stands for the Beer--Lambert law   \cite{w-vang7} when \(t=0\), while stands for the velocity of the fiber pull-back
with increment  \(L=\SI{10}{\milli\metre}\) and time \(t_\mathrm{end}=\SI{10}{\second}\), \textit{i.e.} at a velocity 
\(v = \SI{1.0}{\milli\metre\per\second}\).

\subsection{Heat transfer}

The heat energy was delivered directly to the vein wall \cite{shar}, due to that the compression reduces the vein diameter.
Heat transfer due to the energy of light deposited is described by the following heat and bioheat transfer equations:
\begin{align}\label{flowheat}
\rho c_\mathrm{p} \left( \frac{\partial T}{\partial t} +\mathbf{u}\cdot \nabla T\right)
=\nabla \cdot ( k\nabla T) + q &\qquad \mathrm{in }\ \Omega_\mathrm{lumen}\times ]0 ; t_\mathrm{end}[; \\
\rho c_\mathrm{p}\frac{\partial T}{\partial t} 
 + c_\mathrm{b} \omega (T-T_\mathrm{b}) = \nabla \cdot (k\nabla T) + q  &\qquad \mathrm{in }\  \left(
\Omega_\mathrm{w}\cup \Omega_\mathrm{p}\cup \Omega_\mathrm{s}\right) \times ]0 ; t_\mathrm{end}[,
\label{bioheat}
\end{align}
where:
\begin{description}
\item[ \(T\)] is the temperature [\si{\kelvin}];
\item[ \(\mathbf{u}\)] is the blood velocity vector [\si{\metre\per \second}];
\item[ \(\rho\)] is the density [\si{\kilogram\per\cubic\metre}];
\item[ \( c_\mathrm{p}\)] is the specific heat capacity per unit mass [\si{\joule\per\kilogram\per\kelvin}];
\item[ \(k\)] is the thermal conductivity [\si{\watt\per\metre\per\kelvin}];
\item[ \(q\)] is the heat source caused by laser power [\si{\watt\per\cubic\metre}].
\end{description}
Considering laminar flow, the blood velocity is scalar, and the convective term 
in \eqref{flowheat} reads \(u\frac{\partial T}{\partial z}\). In the heat equation \eqref{flowheat}, 
\(\rho=\rho_\mathrm{b}\)  and \(c_p= c_\mathrm{b}\)  denote the density and the specific heat capacity of the blood, respectively.
The Pennes bioheat transfer equation \eqref{bioheat},   which distinguishes itself  from nonliving systems,
 includes the effects of blood perfusion  \(\omega= \rho_\mathrm{b} w\) [\si{\kilogram\per\cubic\metre\per\second}] 
  that occurs in the capillary bed, that is, the energy transfer term \(- c_\mathrm{b} \omega (T-T_\mathrm{b})\)
is consequence of the mass transport of blood through tissue (cf.  Table \ref{tabthem}).
Here, \(T_\mathrm{b}\) represents the temperature  of the blood (assumed to be \SI{38}{\degreeCelsius}),
 \(w\) denotes  the volumetric  flow [\si{\per\second}],
 and \(c_\mathrm{b} \omega\) accounts  for the heat conducted in direction of the
contribution of flowing blood to the overall energy balance, before the critical coagulation temperature.
\begin{table}   
\fontsize{8pt}{8pt}
\centering
\begin{threeparttable}
\caption{Thermal parameters \cite{marqa,mord}.  }\label{tabthem}   
\small\addtolength{\tabcolsep}{-1pt}
\begin{tabular}{*{27}{c}}
\toprule
\multicolumn{1}{c}{}&
\multicolumn{1}{c}{\textbf{unit}}&
\multicolumn{1}{c}{ \textbf{ blood}}&
\multicolumn{1}{c}{\textbf{vein wall}}&
\multicolumn{1}{c}{\textbf{perivenous tissue}} &
\multicolumn{1}{c}{\textbf{skin}} \\ \midrule
 \(k\) &   \si{\watt\per\metre\per\degreeCelsius}  & 0.52  & 0.53 & 0.21 & 0.21 \\
  \(\rho\) &\si{\kilogram\per\cubic\metre} &    1060  &1080 & 1000 & 1109 \\
\(c_\mathrm{p}\) & \si{\joule\per\kilogram\per\degreeCelsius}  &
 3600 & 3690 & 2350 & 3500  \\
\(\omega\)  & \si{\kilogram\per\cubic\metre\per\second} & & 1.08 & 1 &0.5545  \\
\(A\) &\si{\per\second}& 7.6e+66 & 5.6e+63 & 5.6e+63& 3.1e+98 \\
\(E_\mathrm{a}\) & \si{\joule\per\mol} & 4.48e+05 & 4.30e+05 & 4.30e+05 & 6.28e+05 \\
\bottomrule
\end{tabular}
\end{threeparttable}
\end{table}
 
The heat source \(q\) is induced by the conversion of laser light into heat, the so-called absorbed optical power density,
 since the heat generated by body metabolism is negligible.
The distribution of absorbed energy within the irradiated volume is governed both by the absorption and
the scattering properties of the tissue at the specific wavelength used:
\begin{equation}\label{power}
q=\mu_\mathrm{a} \phi.
\end{equation}
The  heat  exchange on the skin surface is given by the Newton law of cooling
\begin{equation}
k \frac{\partial T}{\partial r}+h_\mathrm{air}( T-T_\mathrm{air})=0
\qquad   \mathrm{if }\  r= r_\mathrm{s}, \ -vt< z<L,\ t>0,
\end{equation}
where \(h_\mathrm{air}\) is the heat transfer coefficient  of the air, and \(T_\mathrm{air}\) denotes the room temperature. In \cite{marqa}, 
\(h_\mathrm{air}\) is assumed to obey an equation that involves the thermal conductivity of the air,
 the characteristic length of flow domain,  and the Prandtl and Reynolds numbers.

On the remaining boundary, no heat transfer outflow is assumed:
\[
\frac{\partial T}{\partial z }=0\qquad \mathrm{if }\ 0<r< r_\mathrm{s}, \ z=-vt,L, \ t>0.
\]

\subsection{Thermal damage to the vein-tissue system}
\label{sdam}

Let \(t_\mathrm{crit}\) be the time correspondent to
  the dimensionless indicator of damage when it is equal to one: \(\Omega(t_\mathrm{crit} )=1\), \textit{i.e.}
  from  the Arrhenius burn integration
\begin{equation}\label{arr}
\frac{1}{A}=\int_{0}^{t_\mathrm{crit}} \exp\left[-\frac{E_\mathrm{a} }{ RT(r,z,\tau)}\right]d\tau,
\end{equation}
where \(R\) is the universal gas constant (\SI{8.314}{\joule\per\mol\per\kelvin}),
\( A\) is a frequency factor [\si{\per\second}], and 
\(E_\mathrm{a}\) is the activation energy for the irreversible damage reaction [\si{\joule\per\mol}]. 
Approximating the above integral by the lower and upper Riemann sums with the partition constituted by a finite number \(M\) of
subintervals of \(] 0, t_\mathrm{crit} [\), of equal length, 
\(t_\mathrm{crit}\) obeys
\[ 
\sum_{m=1}^{M} \exp\left[-\frac{E_\mathrm{a} }{ RT(r,z,m t_\mathrm{crit}/M)}\right] <
\frac{M}{At_\mathrm{crit} }<\sum_{m=1}^{M} \exp\left[-\frac{E_\mathrm{a} }{ RT(r,z,(m-1) t_\mathrm{crit}/M)}\right] .
\] 

Thermal damage is consequence of  the water content of the constituent cells  reaching \SI{ 100}{\degreeCelsius}.
By this reason, the main predictor of the thermal damage is the dense microbubble formation that is commonly seen at the area.
  
\section{Auxiliary solutions}
\label{meth}

There exist several processes of derivation of the required solutions  (see, for instance, \cite{ozis}).
Here,  we firstly use the method of separation of variables to find a family of
elementary solutions to the parabolic equation and the boundary condition, and then the principle of
superposition to construct a solution satisfying the initial condition.
Hereafter, the domain subscripts are dropped out by the sake of simplicity whenever
the meaning of the parameters is well understood in each domain. 

In cylindrical coordinates, the Laplace operator reads  (see, for instance, \cite[p. 9]{ozis})
\[
\Delta=\frac{\partial^2 }{\partial r^2}+ \frac{1}{r}\frac{\partial }{\partial r} +\frac{1}{r^2}
\frac{\partial^2 }{\partial \theta^2}+ \frac{\partial^2 }{ \partial z^2}  ,
\]
with \((r,\theta,z,t)\in ]0;r_\mathrm{o} [\times ]-\pi;\pi[\times  ]0;L[\times ]0;t_\mathrm{end}[\),
where \(r=\sqrt{x^2+y^2}\) and \(\theta\) is the polar angle measured down from the vertical axis \(z\).

Taking the angular symmetry, we seek for solutions of the generic PDE:
\begin{equation}
\label{eqc1cil}
\alpha \frac{\partial \upsilon }{\partial t} + b\frac{\partial \upsilon}{\partial z}- a\left( 
\frac{\partial^2 \upsilon}{ \partial r^2} +
\frac{1}{r}\frac{\partial \upsilon}{\partial r} +
\frac{\partial^2 \upsilon}{ \partial z^2} \right)+ B \upsilon=f, 
\end{equation}
defined in \( ]0;r_\mathrm{o}[\times ]0;L[\times  ]0;t_\mathrm{end}[ \),
with \(\upsilon_0\) denoting the initial datum, \(a >0\), and \(\alpha , b, B \geq 0\).

For our purposes, we begin by exemplifying the decomposition followed for the particular function
\[
f(r,z,t)=f(t) \chi_{[0;r_\mathrm{f} ]} (r) \exp [\iota z ],
\]
where \(\iota\in \mathbb{R}\) and  \(\chi_{[0;r_\mathrm{f} ]} \) stands   for the characteristic function over the interval \([0;r_\mathrm{f} ]\).

Thanks to the Duhamel principle, we look for a  solution which can be of the form
\begin{equation}\label{duhamel}
\upsilon(r,z,t)= v_2(r,z,t)+
\chi_{[0;r_\mathrm{f} ]} (r)  \left(\int_0^t f(s) v_3(r,z,t-s)\mathrm{ds}\right)\exp \left[ \iota z\right] ,
\end{equation}
with \(v_2 = C+v_1\), for some constant \(C\) whenever \(f\equiv BC\),
where \(v_1\) solves (as described in Subsection \ref{analytical}) 
\begin{align}
\alpha\frac{\partial v_1}{\partial t}=a \left( \frac{1}{r}\frac{\partial}{\partial r}\left(r \frac{\partial v_1}{\partial r}\right)+\frac{\partial^2 v_1}{\partial z^2}\right)-
 b\frac{\partial v_1}{\partial z}- B v_1;\label{unsteady} \\
\frac{\partial v_1}{\partial t}(r,-L,t)= 
 \frac{\partial v_1}{\partial z}(r,L,t)=0,\quad\forall r,t;\label{neumL}\\
v_1(r,z,0)=\upsilon_0(r,z), \quad\forall r,z,\label{ic}
\end{align}
and \(v_3\) solves 
\begin{align}\label{v3}
\alpha \frac{\partial v_3}{\partial t} + b\frac{\partial v_3}{\partial z} +
(B+b\iota-a\iota^2) v_3= a\left(
\frac{1}{r}\frac{\partial}{\partial r}\left(r \frac{\partial v_3 }{\partial r}\right)+\frac{\partial^2 v_3}{\partial z^2}\right); \\
v_3 (r,z,0)=1 /\alpha, \quad\forall r,z,\label{v3ci}
\end{align}
such that
\begin{equation}
-a\frac{\partial \upsilon}{\partial r}(r_\mathrm{s},z,t) +h\upsilon(r_\mathrm{s},z,t)=h\gamma,  \quad\forall z,t,
\label{robin}\end{equation}
for some \(\gamma\geq 0\).

We may consider \(v_3(t) =\exp[\zeta t]/\alpha\) with \(\zeta = (a\iota^2-b\iota-B)/\alpha\).

\subsection{Analytical solutions}
\label{analytical}

Using Bernoulli--Fourier technique, the  Cauchy--Robin--Neumann  problem  admits a solution of the form:
\[ 
v_1(r,z,t)= R(r)Z(z) F(t). 
\] 
In order to obtain an analytical solution, let us take the system of ordinary differential equations (ODE)
\begin{equation}\label{Rbeta}
\left\{\begin{array}{lc}
 F'(t) =\zeta F(t)&\\
Z''(z) - (b/a) Z'(z)=\eta^2Z(z)&\\
\left(rR'(r)\right)'=\beta rR(r) ,&\qquad \alpha\zeta= a(\beta+ \eta^2) -B.
\end{array} 
\right.
\end{equation}
The solution of the first ODE is \(F(t)=A_0 \exp[\zeta t]\), for some constant \(A_0\).

The elementary solutions for \(Z\) are
\[
\exp\left[\left(\frac{b } {2a} \pm\Xi \right)z
\right],
\]
where   \(\Xi=(2a)^{-1}\sqrt{b^2+4a^2\eta^2}\).

The elementary solutions for \(R\) are  the Bessel functions of first and second kind and order 0, respectively, \(J_0(\sqrt{|\beta |} r)\) and \(Y_0(\sqrt{|\beta |} r)\)
 if \(\beta <0\);  or the modified Bessel functions of first and second kind and order 0, respectively, 
 \(I_0(\sqrt\beta r)\) and \(K_0(\sqrt\beta r)\) if \(\beta >0\) \cite{ozis}.

We recall their Taylor series expansions around the origin
\begin{align*}
J_0(r) & = \sum_{n=0}^\infty \frac{(-1)^n}{2^{2n}}\frac{r^{2n}}{n!\Gamma(n+1)};\\
Y_0(r) & = \frac{2}{\pi}\left(J_0(r)\ln\left[ \frac{r}{2}\right]-
\sum_{n=0}^\infty (-1)^n\frac{1+(1/2)+\cdots+(1/n)-\gamma}{ 2^{2n}(n!)^2} r^{2n}\right);\\
I_0(r) & = \sum_{n=0}^\infty \frac{1}{2^{2n}}\frac{r^{2n}}{n!\Gamma(n+1)};\\
K_0(r) & =  -I_0(r)\ln\left[ \frac{r}{2}\right]+
\sum_{n=0}^\infty \frac{1+(1/2)+\cdots+(1/n)-\gamma}{2^{2n}(n!)^2} r^{2n},
\end{align*}
where \(\Gamma\) is the gamma function, and \(\gamma\) is
the Euler--Mascheroni constant. Moreover, the following wronskian relationships
\begin{align}
J_1(\beta r)Y_0(\beta r)-Y_1(\beta r)J_0(\beta r) =\frac{2}{\pi\beta r};\label{wronsk1} \\
K_1(\beta r) I_0(\beta r) + I_1(\beta r) K_0(\beta r) =\frac{1}{\beta r} \label{wronsk2}
\end{align}
hold, for any \(\beta>0\) (see \cite[p. 672]{ozis} and \cite[pages 360 and 375]{olver}).

Considering \(\upsilon_0=0\) in \(\Omega_\mathrm{lumen} \) and taking   \( \beta = -\beta_{1,j}^2\) according to  \(j=\)  wall, pad or skin, we find
\begin{align}\label{solf}
v_1(r,z,t) &= 0 
\quad \mathrm{in } \ \Omega_\mathrm{lumen}  ; \\
v_1(r,z,t) &=  A_{1,w}\left( Y_0(\beta_{1,w} r_\mathrm{i}) J_0(\beta_{1,w} r)-J_0(\beta_{1,w} r_\mathrm{i})  Y_0(\beta_{1,w} r) \right)Z(z) \exp[\zeta t] 
\quad \mathrm{in } \ \Omega_\mathrm{w} ; \label{solj}\\
v_1(r,z,t) &=\left( A_{1,p} J_0(\beta_{1,p} r) + A_{2,p} Y_0(\beta_{1,p} r) \right)Z(z) \exp[\zeta t] \quad \mathrm{in } \ \Omega_\mathrm{p} ; \label{soll} \\ \label{sols}
v_1(r,z,t) &=\left( A_{1,s} J_0(\beta_{1,s} r)+ A_{2,s} Y_0(\beta_{1,s} r)  \right)Z(z) \exp[\zeta t]  \quad \mathrm{in } \ \Omega_\mathrm{s}
\end{align}
and  the Neumann condition \eqref{neumL}  implies that
 \begin{equation}\label{defZ}
Z(z)=\exp \left[ \frac{b}{2a } z\right] \left(
b\sinh \left[\Xi (L-z)\right] +\sqrt{b^2+4a^2\eta^2}\cosh \left[\Xi (L-z)\right]\right).
\end{equation}
The above involved constants may be determined by using the continuity conditions on the fluxes and on
the functions themselves.

In particular, by  the homogeneous Robin condition \eqref{robin}, we have
\[
a \beta_{1,s} \left( A_{1,s} J_1(\beta_{1,s} r_\mathrm{s})+ A_{2,s} Y_1(\beta_{1,s} r_\mathrm{s}) \right)
+h\left( A_{1,s} J_0(\beta_{1,s} r_\mathrm{s})+ A_{2,s} Y_0(\beta_{1,s} r_\mathrm{s}) \right) =0.
\]
Notice that if \(\gamma>0\) then we have an additional relation, which can be analyzed by adding an additional solution of the form 
\(R(r)=A_1I_0(\sqrt{B/a}r)+A_2 K_0(\sqrt{B/a}r)\) that obeys
\[
-\sqrt{aB}\left(A_1I_1(\sqrt{\frac{B}{a}}r_\mathrm{s})-A_2 K_1(\sqrt{\frac{B}{a}}r_\mathrm{s})\right)
+h\left(A_1I_0(\sqrt{\frac{B}{a}}r_\mathrm{s})+A_2 K_0(\sqrt{\frac{B}{a}}r_\mathrm{s})\right)=h\gamma.
\]
	
\subsection{Particular solution}
\label{particular}

In this section,  we rephrase the particular solution, according to Duhamel principle,  in \eqref{duhamel} 
\[
 \left(\int_0^t f(s) v_4(r,z,t-s)\mathrm{ds}\right) R(r)\exp \left[ \iota z\right] \]
by taking
\( f(r,z,t)=f(t) R(r) \exp [\iota z ]\),
where \(\iota\in \mathbb{R}\) and \(R\) denotes any of the Bessel functions of order 0, namely, \(J_0(\sqrt{|\beta |} r)\) and \(Y_0(\sqrt{|\beta |} r)\)
 if \(\beta <0\);  or the modified Bessel functions of  order 0, namely,   \(I_0(\sqrt\beta r)\) and \(K_0(\sqrt\beta r)\) if \(\beta >0\),
introduced in Subsection \ref{analytical}.
Analogously, we may consider  a  solution \(v_4(t) =\exp[\zeta t]/\alpha\) with \(\zeta = \left(a(\iota^2+\beta)-b\iota-B\right)/\alpha\).

\section{Results and discussions}

 The primary fluence \(\phi\),  \textit{i.e.} defined in blood whenever the source \(S\not=0\) (cf. \eqref{source}),
can be assumed independent on the radius \(r\) as satisfying \eqref{eqc1cil}  with \(\alpha=1/\nu\),
 \(b=0\), \(a=D\),  \(B=\mu_\mathrm{a}\) and \(f=S\),  with \(\iota = - \mu_\mathrm{t}\). Thus, we may choice the particular functions
\(v_2=0\) and \(v_3 (t) =\nu \exp[\zeta  t]\), 
with
 \[  \zeta:=\nu(D\mu_\mathrm{t} ^2 - \mu_\mathrm{a})> 0
\] according to Table  \ref{tabopt}.
\begin{table}
\fontsize{8pt}{8pt}
\centering
\begin{threeparttable}
\caption{Optical parameters \cite{marqa,ruij}  (at slow shear rate \cite{roggan}).} \label{tabopt}  
\small\addtolength{\tabcolsep}{-1pt}
\begin{tabular}{*{27}{c}}
\toprule
\multicolumn{1}{c}{\textbf{\(\lambda\) [\si{\nano\metre}] }}&
\multicolumn{3}{c}{ \textbf{\(\mu_\textrm{a}\)   [\si{\per\milli\metre} ] }}&
\multicolumn{3}{c}{\textbf{\(\mu_\textrm{s}' \) [\si{\per\milli\metre} ]}}\\
\midrule
\multicolumn{1}{c}{}&
\multicolumn{1}{c}{ blood}& \multicolumn{1}{c}{vein wall}&
\multicolumn{1}{c}{tissue pad} & \multicolumn{1}{c}{skin}&
\multicolumn{1}{c}{ blood}& \multicolumn{1}{c}{vein wall}&
\multicolumn{1}{c}{tissue pad} & \multicolumn{1}{c}{skin} \\ \midrule
810 &0.21& 0.2 & 0.017& 0.2 &  0.73 & 2.4 & 1.2& 0.9\\
980 &0.21 &0.1 & 0.03 &0.10 & 0.6& 2.0& 1.0 & 0.81 \\
1064 & 0.12 & 0.12 & 0.034 & 0.10 & 0.58 & 1.95 & 0.98 & 0.77\\
\bottomrule
\end{tabular}
\end{threeparttable}
\end{table}
Then,  we have 
\[ 
\phi_\nu(r,z,t) = \frac{\nu S(r_\mathrm{f}, z, 0)}{\zeta +\mu_\mathrm{t} v }
\exp[\zeta t]\left(1- \exp [ - (\zeta + \mu_\mathrm{t} v) t]\right)
\] 
for \(  0\leq r < r_\mathrm{f}\), \(-L< z < L\),  \(  0\leq t < t_\mathrm{end}\).
Typically, a  laser fiber for medical applications has a 600-micron diameter (\(r_\mathrm{f}=\SI{0.3}{\milli\metre}\)).

Next,  by the interface continuity conditions,  \(\phi_\nu\) can be extended as a solution
at the position \((r,z)\) and the time \(t\),
\begin{align*}
\phi_\nu (r,z,t)  & =  \frac{\nu S(r_\mathrm{f}, z, 0 )}{\zeta + \mu_\mathrm{t} v}\exp [\zeta t] \\
&+\big(  B_1J_0(\beta_j r) + B_2 Y_0(\beta_j r) \big) \exp [ -\mu_\mathrm{t} (z + vt) ]
 \quad\mathrm{if }\ r_\mathrm{f}< r \leq r_\mathrm{i} \ (j= \mathrm{blood});\\
\phi_\nu (r,z,t) & = \exp [ -\mu_\mathrm{t}z ] \Big( 
 ( B_3 I_0(\varkappa_j r) + B_4 K_0(\varkappa_j r) ) \exp [ \zeta t] \\
&+ (B_5 J_0(\beta_j r) + B_6 Y_0(\beta_j r) ) \exp [- \mu_\mathrm{t} vt] \Big)
\quad\mathrm{otherwise},
\end{align*}
where  the abstract constants \(B_1,\cdots,B_6\) are defined  by the boundary and interface continuity conditions.The parameters,
\(\varkappa_j\) and \(\beta_j\),  are
determined due to that the PDE \eqref{light} is verified by the function \(\phi\), namely
\begin{align*}
\zeta/\nu_j + \mu_\mathrm{a} ^{(j)} =  D_j \left(\mu_\mathrm{t}  ^2 + \varkappa_j^2\right)\quad
j= \mathrm{wall}, \mathrm{pad}, \mathrm{skin} ;\\
-\frac{\mu_\mathrm{t}v}{\nu_j} + \mu_\mathrm{a} ^{(j)} =  D_j \left(\mu_\mathrm{t}  ^2 - \beta_j^2\right)\quad
j= \mathrm{blood}, \mathrm{wall}, \mathrm{pad}, \mathrm{skin}.
\end{align*}
That is,
\begin{align*}
\varkappa_j &= \sqrt{
 \frac{ n_\mathrm{tissue} }{n_\mathrm{blood}}\frac {D\mu_\mathrm{t} ^2 - \mu_\mathrm{a} }{ D_j }+ (\mu^{( j )}_\mathrm{eff})^2 - \mu_\mathrm{t} ^2} >0 \qquad 
 j=  \mathrm{wall}, \mathrm{pad}, \mathrm{skin} ;\\
\beta_j &=\sqrt{ \mu_\mathrm{t} ^2  -(\mu^{( j )}_\mathrm{eff})^2 + \frac{  \mu_\mathrm{t} v }{ \nu_j D_j }} >0 \qquad
j= \mathrm{blood}, \mathrm{wall}, \mathrm{pad}, \mathrm{skin}.
\end{align*}
Here, we may consider \(n=1.4\) for both the blood and the tissues.

This solution proves that the problem is ill-posed.
Although it gives a good answer at the scale of picosecond   (\(c=\SI{0.3}{\milli \metre\per \pico\second}\)), 
it is  inadequate for describing the behavior of the fluence rate whenever the fiber moves.
The unsteady \(\phi\) should solves \eqref{light} at the steady state, being such that \(\phi\) attains its maximum at \(z=-vt\).
Then, we have 
\begin{align}\label{phi1}
\phi(r,z,t)  & =
B_0\exp[-\mu_\mathrm{eff}  (z+vt)]-\frac{S(r_\mathrm{f}, 0,0 )}{D\mu_\mathrm{t} ^2 - \mu_\mathrm{a}} \exp [- \mu_\mathrm{t}(z + vt) ]\quad \mathrm{if }\  
0\leq r < r_\mathrm{f};\\
\phi(r,z,t)  & =  B_0\exp [-\mu_\mathrm{eff} (z+vt)] \nonumber\\
&+\big(  B_1J_0(\beta_\mathrm{b} r) + B_2 Y_0(\beta_\mathrm{b} r) \big) \exp [ -\mu_\mathrm{t} (z + vt) ]
 \quad\mathrm{if }\ r_\mathrm{f}< r \leq r_\mathrm{i} ;\\
\phi(r,z,t) & =  ( B_3 W_3(\varkappa_{\lambda,j} r) + B_4W_4(\varkappa_{\lambda,j} r) )  \exp [ -\mu_\mathrm{eff}(z+vt) ] \nonumber\\
& + (B_5 J_0(\beta_j r) + B_6 Y_0(\beta_j r) ) \exp [- \mu_\mathrm{t} (z + vt )] 
\quad\mathrm{otherwise},\label{phi3}
\end{align}
where  the abstract constants \(B_0,B_1,\cdots,B_6\) are defined  by the initial,  boundary and interface continuity conditions.
We consider the modified Bessel functions  \(W_3 =I_0\) and \(W_4 =K_0\) if \(\lambda = 810;1064\) and \(j=\) wall, skin, or \(\lambda = 980\) and \(j= \) wall;
and the Bessel functions  \(W_3 =J_0\) and \(W_4 =Y_0\) otherwise, in accordance to the factors of the \(r\)-argument
\begin{align}
\varkappa_{\lambda,j} &=\left\{\begin{array}{ll}
 \sqrt{ (\mu^{( j )}_\mathrm{eff})^2 - \mu_\mathrm{eff} ^2} >0 & j=  \mathrm{wall}, \mathrm{skin} ;\\
 \sqrt{  \mu_\mathrm{eff} ^2 - (\mu^{( j )}_\mathrm{eff})^2} >0 & j=\mathrm{pad},
\end{array}
\right.  \mbox{ if }\lambda = 810;1064\\
\varkappa_{\lambda,j} &=\left\{\begin{array}{ll}
 \sqrt{ (\mu^{( j )}_\mathrm{eff})^2 - \mu_\mathrm{eff} ^2} >0 & j=  \mathrm{wall};\\
 \sqrt{  \mu_\mathrm{eff} ^2 - (\mu^{( j )}_\mathrm{eff})^2} >0 & j=\mathrm{pad}, \mathrm{skin} ,
\end{array}
\right.  \mbox{ if }\lambda = 980\\
\beta_j &=\sqrt{ \mu_\mathrm{t} ^2  -(\mu^{( j )}_\mathrm{eff})^2} >0 \qquad
j= \mathrm{blood}, \mathrm{wall}, \mathrm{pad}, \mathrm{skin}.\label{beta}
\end{align}

Parameters, used in clinical procedures, are known:
the power is set at \SI{15}{\watt} (with wavelengths of \SI{810}{\nano\metre} and \SI{980}{\nano\metre}) 
and at \SI{10}{\watt} (with wavelengths of \SI{980}{\nano\metre} and \SI{1064}{\nano\metre}). 
\begin{figure}\begin{multicols}{2}
\centering
\includegraphics[width=0.5\textwidth]{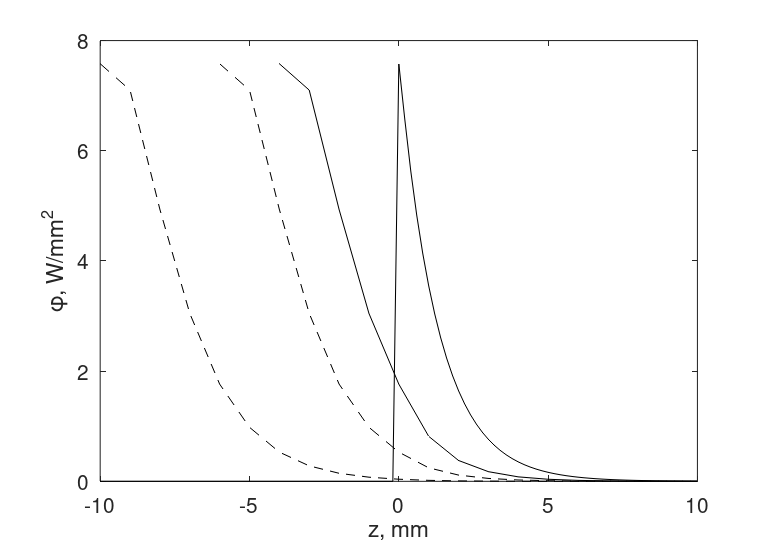}\\
\includegraphics[width=0.5\textwidth]{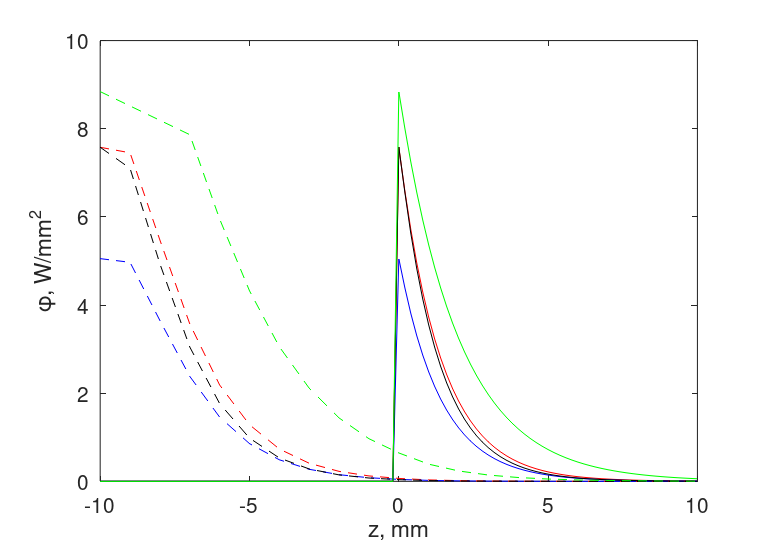}
\end{multicols}
\caption{\footnotesize Left: Graphical representations of \(\phi\), at the  wavelength of \SI{810}{\nano\metre}, for different instants of time.
Right: Graphical representations of \(\phi\) for the wattage set at \SI{15}{\watt}:
with wavelengths of \SI{810}{\nano\metre} (in black) and \SI{980}{\nano\metre} (in red)
and at  \SI{10}{\watt}: with wavelengths of \SI{980}{\nano\metre} (in blue) and \SI{1064}{\nano\metre} (in green).
Solid lines stand for the initial instant of time, while dashed lines stand for the final instant of time.}
\vspace{0.2cm}
\hrule
\label{fsource}
\end{figure} 
Calculations use Octave software, under the optical parameters  in Table \ref{tabopt}.
Figure \ref{fsource} (Left) shows the \(\phi\)-profile to the power of \SI{15}{\watt} and the wavelength of  \SI{810}{\nano\metre} at different instants of time,
considering the continuous movement of the fiber tip.
The slope decreases from  the initial instant of time (solid line) \(z=t=0\)  until the final instant of time (dashed line) \(z=-10, t=10\),
which reflects the accumulation of the fluence rate.
Figure \ref{fsource} also illustrates the higher distributions  only under the tip, as it is experimentally consistent.

As expected,  Figure \ref{fsource} (Right) shows similar profiles between different wavelengths and set powers.
At the wavelength  of  \SI{980}{\nano\metre},  the \SI{15}{\watt}-curves (in red) has higher values than the  \SI{10}{\watt}-curves (in blue).
At the power  \SI{15}{\watt}, the  \SI{810}{\nano\metre}-curves (in black) and the  \SI{980}{\nano\metre}-curves (in red) match each other,
while at the power  \SI{10}{\watt},  the slope of  \SI{980}{\nano\metre}-curves (in blue) is similar than the slope of 
 \SI{1064}{\nano\metre}-curves (in green), but there is no match.

This evaluation improves the study of the absorbed volumetric power \(q\), and consequently it
will improve the study of the distribution of the temperature and make the simulation of  the  heat transfer closer to reality.
Our results are consistent to that  the use of laser energies of various wavelengths has
no significant difference in their effectiveness and complication rate \cite{meme}.

Hereafter, we assume the solution \(\phi\) as defined by \eqref{phi1}-\eqref{beta}.
According to Section \ref{meth},
let \(T -T_\mathrm{b}\) be the thermal  solution of \eqref{eqc1cil} with \(\alpha = \rho c_\mathrm{p}\),  \(a=k\), \(f =q\), and
\begin{itemize}
\item in the domain \(\Omega_\mathrm{lumen}\): \(b = \rho_\mathrm{b} c_\mathrm{b} u\) and \(B=0\);
\item in the domain \(\Omega_\mathrm{w}\cup \Omega_\mathrm{p}\cup\Omega_\mathrm{s}\):  \(b=0\)  and \(B= c_\mathrm{b} \omega\).
\end{itemize}

Firstly,  at the position \((r,z)\in [0;r_\mathrm{f}[\times ]-L;L[\) and the time \(t>0\),  we split into
\begin{description}
\item[(i)]  \(f(t) = \mu_\mathrm{a} \frac{S(r_\mathrm{f}, 0,0 )}{D\mu_\mathrm{t} ^2 - \mu_\mathrm{a}}\exp[-\mu_\mathrm{eff}vt]\) and \(\iota =-\mu_\mathrm{eff}\);
\item[(ii)]  \(f(t) =- \mu_\mathrm{a} \frac{S(r_\mathrm{f}, 0,0 )}{D\mu_\mathrm{t} ^2 - \mu_\mathrm{a}}\exp[-\mu_\mathrm{t}vt]\) and \(\iota =-\mu_\mathrm{t}\).
\end{description}
Analogously in the determination of \(\phi_\nu\) we use \(v_3(t)= \exp[\zeta t]/(\rho_\mathrm{b}c_\mathrm{b})\) solving \eqref{v3}-\eqref{v3ci} with
the factors of the \(t\)- argument
\begin{description}
\item[(i)] \(\zeta_1 = k\mu_\mathrm{eff}^2/(\rho_\mathrm{b}c_\mathrm{b})+u\mu_\mathrm{eff}\).
\item[(ii)] \(\zeta _2= k\mu_\mathrm{t}^2/(\rho_\mathrm{b}c_\mathrm{b})+u\mu_\mathrm{t}\).
\end{description}

Secondly,  for \(r_\mathrm{f}< r<r_\mathrm{i}\), we use  \(v_4(t)= \exp[\zeta_3 t]/(\rho_\mathrm{b}c_\mathrm{b})\) from Subsection \ref{particular} with
 \(\zeta _3= k\mu_\mathrm{eff}^2/(\rho_\mathrm{b}c_\mathrm{b})+u\mu_\mathrm{t}\) by
taking \(\beta_\mathrm{b}\) that is given in \eqref{beta} into account. 
Next, we similarly argue for the domain \(\Omega_\mathrm{w}\cup \Omega_\mathrm{p}\cup\Omega_\mathrm{s}\),
 concluding
\begin{align}\label{Trf}
T(r,z,t) &= T_\mathrm{b}+  
\mu_\mathrm{a}B_0\exp[-\mu_\mathrm{eff} z]
\frac{ \exp [\zeta_1 t]-\exp[-\mu_\mathrm{eff}vt]  }{k\mu_\mathrm{eff}^2+\rho_\mathrm{b}c_\mathrm{b}  \mu_\mathrm{eff}(u+v) } \nonumber\\
& -\frac{ \mu_\mathrm{a} S(r_\mathrm{f}, 0,0) }{D\mu_\mathrm{t} ^2 - \mu_\mathrm{a} }\exp[-\mu_\mathrm{t} z]  \frac{ \exp [\zeta_2 t]-\exp[- \mu_\mathrm{t} vt] }{
k\mu_\mathrm{t} ^2 + \rho_\mathrm{b}c_\mathrm{b}  \mu_\mathrm{t}(u+v)}   
\quad \mathrm{ if }\  0\leq r\leq r_\mathrm{f};\\
	T(r,z,t) &=  T_\mathrm{b}+ \mu_\mathrm{a}  B_0
\frac{ \exp [-\mu_\mathrm{eff}( z-( k\mu_\mathrm{eff}/(\rho_\mathrm{b}c_\mathrm{b})+u)t)]
-\exp[-\mu_\mathrm{eff} (z+vt)]  }{k\mu_\mathrm{eff}^2+\rho_\mathrm{b}c_\mathrm{b}  \mu_\mathrm{eff}(u+v) } \nonumber\\
&+\big(  B_1J_0(\beta_\mathrm{b} r) + B_2 Y_0(\beta_\mathrm{b} r) \big) \exp[ - \mu_\mathrm{t} z ]
\frac{\exp [ \zeta_3 t]-\exp[- \mu_\mathrm{t} vt]  }{k\mu_\mathrm{eff} ^2  + \rho_\mathrm{b}c_\mathrm{b}  \mu_\mathrm{t}(u+v) } \nonumber \\ 
&+T_1(r,z,t) \qquad\mathrm{if }\ r_\mathrm{f}< r \leq r_\mathrm{i}; \\
	T(r,z,t) &=  T_\mathrm{b} +( B_3 W_3(\varkappa_{\lambda,j} r) + B_4W_4(\varkappa_{\lambda,j} r) )  \exp [ -\mu_\mathrm{eff}z ] \nonumber\\
&\qquad\times\frac{ \exp \left[\sqrt{(k(\mu_\mathrm{eff}^ {(j)})^2-\rho_\mathrm{b}c_\mathrm{b}w)/(\rho c_\mathrm{p})} t\right]-\exp[-\mu_\mathrm{eff}vt]
 }{k(\mu_\mathrm{eff}^ {(j)})^2- c_\mathrm{b} \omega } \nonumber\\
& + (B_5 J_0(\beta_j r) + B_6 Y_0(\beta_j r) ) \exp [- \mu_\mathrm{t} z] \nonumber\\
&\qquad\times
\frac{ \exp \left[\sqrt{(k(\mu_\mathrm{eff}^ {(j)})^2-\rho_\mathrm{b}c_\mathrm{b}w)/(\rho c_\mathrm{p})} t\right]-\exp[- \mu_\mathrm{t} vt] 
 }{k(\mu_\mathrm{eff}^ {(j)})^2- c_\mathrm{b} \omega } \nonumber\\
&+T_1(r,z,t) \qquad\mathrm{otherwise},\label{Trs}
\end{align}
where \(T_1\) is the combination of radial dependent Bessel functions and longitudinal and temporal dependent exponential functions
such that   \eqref{unsteady}-\eqref{ic} as well as the interface continuity conditions are  verified.

As the lumen \(\Omega_\mathrm{lumen}\) is constituted by the blood, two different situations exist: 
\begin{description}
\item[Case 1.] The blood flow is obstructed (\(u=0\)), for instance the vein is completely clamped or the SSV in  the presence of the inserted catheter.
 \item[Case 2.] The blood flows at \(u=\SI{70}{\milli\metre\per\second}\), as such it happens in the GSV. 
The diameter of the GSV varies from \SIrange{11}{12}{\milli\metre} at SFJ until  \SIrange{7.5}{8}{\milli\metre} at the proximal thigh
(at the knee level)  \cite{doganci}.
\end{description}

In the case 1, \(T_1\) is given by \eqref{solf}-\eqref{sols} with \(Z(z) =\cosh[\eta (L-z)]\) and
\[
\rho_j c_{\mathrm{p},j}\zeta = - k_j (\beta_{1,j}^2 +\eta^2) - c_{\mathrm{p} ,j}\omega_j \qquad j =\mathrm{w},\mathrm{p},\mathrm{s}.
\]

In the case 2, the solution \(T_1\) depends on the general form  \eqref{defZ} of \(Z\), and \eqref{solf}-\eqref{sols}
 take the dependence on the vein size into account. Known values exist
for diameters ranging between 3 and \SI{10}{\milli\metre}, or greater than \SI{10}{\milli\metre}, according to
the small saphenous vein  (SSV), anterior accessory vein, and great  saphenous vein,
namely, diameters of \SI{7.5}{\milli\metre} \cite{goode}, \(r_\mathit{i}=\SI{3.75}{\milli\metre}\)
 and \(\varepsilon =\SI{0.75}{\milli\metre}\), 
and of \SI{1.2}{\centi\metre} \cite{gold}, \(r_\mathit{i}=\SI{6}{\milli\metre}\) and \(\varepsilon =\SI{1.2}{\milli\metre}\).
The proposed solution may address the quantitative questions
in the context of the thermal ablation  treatment under study.
For    elucidating the effect of vein diameter, further study will be the aim of future work.

Finally, we may use the formula \eqref{Trf}-\eqref{Trs} to read the temporal course of damage events from the spatial domain.
However, a first analysis should be done.
Considering in \eqref{arr} that \(T\geq T_\mathrm{min}\), we find the following upper bound 
\[
t_\mathrm{crit}\leq \frac{1}{A}\exp\left[\frac{E_\mathrm{a}}{RT_\mathrm{min}}\right].
\]
Next, if we use the thermal parameters from Table \ref{tabthem}, 
the above upper bound of  the critical time \(t_\mathrm{crit}\) can be calculated  (cf. Table  \ref{tabtemperature})
 in function of different the minimum surface temperatures.
\begin{table}    
\fontsize{8pt}{8pt}
\centering
\begin{threeparttable}
\caption{Upper bounds of \(t_\mathrm{crit}\) in seconds. }  \label{tabtemperature}  
\small\addtolength{\tabcolsep}{-1pt}
\begin{tabular}{*{27}{c}}
\toprule
\multicolumn{1}{c}{ \si{\degreeCelsius} }&
\multicolumn{1}{c}{ \textbf{ blood}}&
\multicolumn{1}{c}{\textbf{vein wall}}&
\multicolumn{1}{c}{\textbf{perivenous tissue}} &
\multicolumn{1}{c}{\textbf{skin}} \\ \midrule
50& 3.4e+05  & 5.8e+05  & 5.8e+05  & 1.1e+03\\
60&   2.3e+03 &  4.7e+03 &  4.7e+03  & 9.5e-01\\
70&   2.1e+01 &  5.1e+01 &  5.1e+01 &  1.3e-03\\
 80&  2.4e-01 &  7.2e-01 &  7.2e-01 &  2.5e-06\\
 90&  3.6e-03  & 1.3e-02  & 1.3e-02 &  6.9e-09\\
 100&  6.8e-05  & 2.8e-04 &  2.8e-04  & 2.6e-11\\
\bottomrule
\end{tabular}
\end{threeparttable}
\end{table}
Then, the operating time of \SI{10}{\second} is a safe value for the thermal thresholds for tissue damage 
at the vein wall, perivenous tissue  or skin., the so-called damage temperature (\(T_\mathrm{min}= \SI{50}{\degreeCelsius}\)).
This value is consensual among clinicians and researchers  in the ablation treatments,  \cite{cs} and the references therein.
We conclude that the  interpretation of the vein-tissue system damage substantially disagrees with the critical temperature of \SI{50}{\degreeCelsius}
being the  temperature that temperatures above it result in necrosis.

The present result shows that the blood coagulates before the dehydration/necrosis of the wall tissue. 
Moreover,  for the blood threshold \(T_\mathrm{min}= \SI{100}{\degreeCelsius}\),  the operating time clearly surpasses the upper bound, which is consistent with 
that  a thin layer of carbonized blood  is found to cover the fiber tip.
Although according to  \cite{ruij} 
the black layer occurs at temperatures around  \SI{300}{\degreeCelsius}.  This black layer
absorbs an average of 45\% of the emitted light power, resulting in a decrease of the tip temperatures.

Indeed, much work remains to be done.
The damage  of the vein reduces to the vaporisation and occlusion in both above situations, namely cases 1 and 2.
The failure and complication rates depend on vein size  \cite{chaar}.
Also, the perivenous tumescence injection (tumescence anesthesia) is carried on  to protect
the perivenous tissue from thermal damage and reduce the lumen of the truncal vein by compression and spasm \cite{schm}.

\section{Conclusions}

The  derived solutions, namely \(\phi\) and \(T\),  might be a tool to generate quantitative and/or qualitative results.
Besides, the configurations of the physical problem solved are questioned.
Our main conclusion is that the parabolic  equation for the light transport leads to the application of the pulsed laser of the order of picoseconds only,
 while the elliptic equation leads to the physical solution for a continuous laser light.

\section*{Acknowledgements}
 Deeply thanks to Professor Lu\'\i s Filipe V. Ferreira by awakening my interest on the light propagation field.

\end{document}